\documentclass{article}

\title{Robust finite element error estimates for elliptic and parabolic
distributed optimal control problems \\ with energy regularization} 
\author{Ulrich~Langer\footnote{Johann Radon Institute for Computational and
    Applied Mathematics, Austrian Academy of Sciences, Altenberger Stra{\ss}e 69, 4040 Linz, Austria, Email:
    ulrich.langer@ricam.oeaw.ac.at}, 
    \; Olaf~Steinbach\footnote{Institut f\"{u}r Angewandte Mathematik,
    Technische Universit\"{a}t Graz, Steyrerga{\ss}e 30, 8010 Graz, Austria,
    Email: o.steinbach@tugraz.at}, 
    \; Huidong~Yang\footnote{Johann Radon Institute for Computational and Applied
    Mathematics, Austrian Academy of Sciences, Altenberger Stra{\ss}e 69, 4040 Linz, Austria, Email:
    huidong.yang@ricam.oeaw.ac.at} 
}  

\date{\today}

\usepackage{amsmath}
\usepackage{amsthm}
\usepackage{amsfonts}
\usepackage{booktabs}
\usepackage{graphicx}
\usepackage{diagbox}

\newtheorem{lemma}{Lemma}

\begin{document}
\title{Space-time finite element methods for the initial temperature
  reconstruction}
\author{Ulrich~Langer\footnote{Johann Radon Institute for Computational and
    Applied Mathematics, Austrian Academy of Sciences,
    Altenberger Stra{\ss}e 69, 4040 Linz, Austria, Email:
    ulrich.langer@ricam.oeaw.ac.at}, 
    \; Olaf~Steinbach\footnote{Institut f\"{u}r Angewandte Mathematik,
    Technische Universit\"{a}t Graz, Steyrergasse 30, 8010 Graz, Austria,
    Email: o.steinbach@tugraz.at}, 
  \; Fredi~Tr\"oltzsch\footnote{Institut f\"ur Mathematik,
    Technische Universit\"at Berlin, Stra{\ss}e des 17.~Juni 136,
    10623 Berlin, Germany, Email: troeltzsch@math.tu-berlin.de},
    \; Huidong~Yang\footnote{Johann Radon Institute for Computational
    and Applied
    Mathematics, Austrian Academy of Sciences, Altenberger
    Stra{\ss}e 69, 4040 Linz, Austria, Email:
    huidong.yang@ricam.oeaw.ac.at} 
}  

%
%
\maketitle

\abstract{This work is devoted to the  reconstruction of the initial temperature
  in the backward heat equation using the space-time finite element method on fully
  unstructured space-time simplicial meshes proposed by Steinbach (2015). Such
  a severely ill-posed problem is tackled by the standard Tikhonov
  regularization method. This leads to a related optimal control for an parabolic equation
  in the space-time domain. In this setting, the control is taken as
  initial condition, whereas the terminal observation data 
  serve as target. 
  The objective becomes a standard terminal observation functional
  combined with the Tikhonov regularization.
  The space-time finite element method is applied to the space-time optimality 
  system that is well-posed for a fixed regularization parameter.
}
%
%
\section{Introduction}\label{intro}
In this work, we investigate the applicability of unstructured space-time
methods to the numerical solution of inverse problems using the classical
inverse problem of the reconstruction of the initial temperature in the heat
equation from an observation of the temperature at a finite time horizon: Find
the initial temperature $u_0(\cdot) := u(\cdot,0) \in L^2(\Omega)$  on
$\Sigma_0$ such that 
\begin{equation}\label{backheat2Dstrong}
  \partial_t u - \Delta_x u = 0 \quad \mbox{in} \; Q,\quad
  u=0 \quad \mbox{on} \; \Sigma, \quad
  u=u_T^\delta \quad \mbox{on}\; \Sigma_T \, ,
\end{equation}
where $Q := \Omega \times (0, T)$ denotes the space-time cylinder with the
boundary $\partial Q = \overline{\Sigma} \cup \overline{\Sigma}_0 \cup
\overline{\Sigma}_T$, $\Sigma:= \partial\Omega\times(0, T)$, $\Sigma_0
:=\Omega \times \{0\}$, $\Sigma_T:= \Omega \times {T}$, the bounded Lipschitz
domain $\Omega \subset {\mathbb{R}}^d$, $d \in \{1,2,3\}$, and a finite time
horizon $T > 0$. Moreover, $u_T^\delta  \in L^2(\Omega)$ denotes the observed
terminal temperature which may contain some noise characterized by the noise level $\delta \ge 0$,
\begin{equation}\label{eq:noiselevel}
  \| u_{T}^\delta - u_{T} \|_{L^2(\Omega)} \leq \delta, 
\end{equation}
where $u_{T} = u(\cdot,T) \in L^2(\Omega)$ represents the unpolluted exact
data.  

In contrast to the forward heat equation with known initial data, the backward
heat equation (\ref{backheat2Dstrong}) is severely ill-posed; see
\cite[Example~2.9]{EnglHankeNeubauer1996}. In fact, the solution of
(\ref{backheat2Dstrong}) does not continuously depend on the data
$u_T^{\delta}$ even when the solution exists. Following the notation in
\cite{EnglHankeNeubauer1996}, the problem (\ref{backheat2Dstrong})
may be reformulated as an abstract operator equation in a more general
setting: Find $u_0\in{\mathcal X}$ such that
\begin{equation}\label{opteqn}  
  S u_0=u_T,  
\end{equation}
where $S:{\mathcal X}\rightarrow{\mathcal Y}$ denotes a bounded linear
operator between two Hilbert spaces ${\mathcal X}$ and ${\mathcal Y}$. It is
clear that there does not exist a continuous inverse operator
$S^{-1}:{\mathcal  Y}\rightarrow{\mathcal X}$ in general. Therefore, we
consider a regularized solution, depending on the choice of Tikhonov's
regularization parameter $\varrho:=\varrho(\delta)$,  
\begin{equation*}
  u_{0}^{\delta, \varrho}:=\left(S^{*}S+\varrho I\right)^{-1}S^{*}u_{T}^\delta ,
\end{equation*}
as the unique minimizer of the Tikhonov functional
\cite{TikhonovArsenin:1977} 
\begin{equation}\label{redopt}
  {\mathcal J}_\varrho(z) := \frac{1}{2} \, \|Sz-u_T^\delta\|_{\mathcal Y} +  
  \frac{\varrho}{2} \, \|z \|^2_{\mathcal X}.
\end{equation}
It is well known that we have  the  convergence
\[
\lim_{\delta\rightarrow 0} u_{0}^{\delta, \varrho} = u_{0}^\dagger
    {\text{ in }{\mathcal X},}   \text{ if the conditions} 
    \;\;
    \lim\limits_{\delta\rightarrow 0} \varrho(\delta)=0 
    \;\; \mbox{and} \;
    \lim_{\delta\rightarrow 0}\frac{\delta^2}{\varrho(\delta)} = 0
\]
are satisfied. Here, $u_0^\dagger$ denotes the best-approximated solution to
the operator equation (\ref{opteqn}); see \cite[Theorem
  5.2]{EnglHankeNeubauer1996} for a more detailed discussion, and also
\cite{ChenHofmannZou2018,LeykekhmanVexlerWalter2020}. 

The main focus of this work is to describe a space-time finite element method
(FEM) on fully unstructured simplicial meshes to solve the minimization
problem (\ref{redopt}) subject to the solution of the heat equation
(\ref{backheat2Dstrong}). Such a space-time method has been studied for
the forward heat equation in \cite{OS15}, and for other parabolic optimal
control problems in
\cite{LangerOlafTroltzschYang01,LangerOlafTroltzschYang}.  

The remainder of this paper is structured as follows: In
Section~\ref{modelprob}, we discuss the related optimal control
problem. Its solution is obtained by the optimality system consisting of the
(forward) heat equation, the adjoint heat equation, and the gradient
equation. Based on the Banach--Ne\v{c}as--Babu\v{s}ka theory
\cite{ErnGuermond2004}, we establish unique solvability of the resulting
coupled system, when eliminating the unknown initial datum. In
Section~\ref{dist}, for the numerical solution of the inverse problem
\eqref{backheat2Dstrong}, we first consider the discrete optimal control
problem, which is based on the space-time discretization of the forward
problem. The solution is characterized by a discrete gradient equation, which
turns out to be the Schur complement system of the discretized coupled
variational formulation. First numerical results  are reported in Section
\ref{numa}. These results show the potential of the space-time approach
proposed. Finally, some conclusions are drawn in Section~\ref{sec:con}. 

%
%
\section{The related optimal control problem}
\label{modelprob}
In our case, the Hilbert spaces ${\mathcal X}$ and ${\mathcal Y}$ are
specified as ${\mathcal X} = {\mathcal Y} = L^2(\Omega)$, and the image $Sz$
of the operator $S:L^2(\Omega)\to L^2(\Omega)$  in the Tikhonov functional
\eqref{redopt} is defined by the solution $ u \in X:= L^2(0,T;
H_0^1(\Omega)) \cap H^1(0,T;H^{-1}(\Omega))$ of the forward heat conduction
problem 
\begin{equation}\label{conheat}
  \partial_t u - \Delta_x u = 0 \quad \mbox{in} \; Q, \quad
  u=0 \quad \mbox{on} \; \Sigma,\quad u=z \quad \mbox{on}\; \Sigma_0,
\end{equation}
and its evaluation on $\Sigma_T$, i.e., $(Sz)(x) = u(x,T)$,
$x\in\Omega$. Here, the control $z \in L^2(\Omega)$ represents the initial
data in \eqref{conheat}. Rewriting the minimization of the functional
\eqref{redopt} in terms of $z$, we obtain the {\em optimal control
  problem}
\begin{equation}\label{fullopt}
  {\mathcal J}_\varrho(z) 
  := \frac{1}{2} \,
  \| u(x, T)-u_T^\delta\|^2_{L^2(\Omega)} +
  \frac{\varrho}{2} \, \|z\|^2_{L^2(\Omega)} \to
  \min_{z \in L^2(\Omega)},
\end{equation}
where the {\em state} $u \in X$ is associated to the {\em control} $z$ subject
to \eqref{conheat}.   

To set up the necessary and sufficient optimality conditions for the optimal
control $z$ with associated state $u$, we introduce the adjoint equation 
\begin{equation}\label{adjoint}
  -\partial_t p - \Delta_x p = 0 \quad \mbox{in} \; Q, \quad
  p=0 \quad \mbox{on} \; \Sigma, \quad
  p = u - u_T^\delta \quad \mbox{on} \; \Sigma_T .
\end{equation}
It has a unique solution $p \in X$, the {\em adjoint state}. The adjoint
equation can be derived by a formal Lagrangian technique as in
\cite{Troltzsch2010}. If $z$ is the optimal control with associated state
$u\in X$, then a unique adjoint state $p\in X$ solving \eqref{adjoint}
exists such that the {\em gradient equation}
\begin{equation}\label{kktz}
  p + \varrho \, z = 0 \quad \textup{on } \Sigma_0
\end{equation}
is satisfied. Using this equation, we can eliminate the unknown initial datum
$z$ in the state equation \eqref{conheat} to conclude  
\begin{equation}\label{conheat 1}
  \partial_t u - \Delta_x u = 0 \quad \mbox{in} \; Q, \quad
  u=0 \quad \mbox{on} \; \Sigma,\quad u= - \frac{1}{\varrho} \, p
  \quad \mbox{on}\; \Sigma_0
\end{equation}
for the optimal state $u$. The {\em reduced optimality system}
\eqref{adjoint},\eqref{conheat 1} is necessary and sufficient for
optimality of $u$ with associated adjoint state $p$. In what follows, we will
describe a space-time finite element approximation of this system.

The space-time variational formulation of the heat equation in
\eqref{conheat 1} (without initial condition) is to find $u \in X$ such
that 
\begin{equation}\label{Def b}
  b(u,v) := \int_0^T \int_\Omega \Big[
  \partial_tu(x,t) v(x,t)+\nabla_xu(x,t)\cdot\nabla_x v(x,t)
  \Big] \, dx \, dt = 0
\end{equation}
is satisfied for all $v \in Y:=L^2(0,T;H_0^1(\Omega))$. The spaces $X$ and $Y$
are equipped with the norms 
 \[
  \| v \|_Y = \| \nabla_x v \|_{L^2(Q)}
  \;\; \mbox{and} \;\;
  \| u \|_X = \sqrt{\| \partial_t u \|_{Y^*}^2 + \| u \|_Y^2} =
  \sqrt{\| w_u \|_Y^2 + \| u \|_Y^2},
\]
with $w_u \in Y$ being the unique solution of the variational problem
\[
  \int_0^T \int_\Omega \nabla_x w_u(x,t) \cdot \nabla_x v(x,t) \, dx \, dt
  = \int_0^T \int_\Omega \partial_t u(x,t) \, v(x,t) \, dx \, dt \quad
  \forall\, v \in Y .
\]
We multiply the adjoint heat equation \eqref{adjoint} by a test function
$q \in X$, integrate over $Q$, and apply integration by parts both in space
and time. Then we insert the  terminal data $u(T) - u_T^\delta$ of $p$  in the
arising term $p(T)$, and substitute the term $p(0)$ by $- \rho z = - \rho
u(0)$ in view of \eqref{kktz}. In this way, we arrive at the weak form of
the adjoint problem \eqref{adjoint} 
\begin{eqnarray*}
  0
  & = & \int_0^T \int_\Omega \Big[ - \partial_t p(x,t) \, q(x,t) -
        \Delta_x p(x,t) \, q(x,t) \Big] dx \, dt \\
  & = & - \int_\Omega [u(x,T)-u^\delta_T(x)] \,  q(x,T) \, dx
        - \varrho \int_\Omega u(x,0) \, q(x,0) \, dx\\
  & & \hspace*{1cm} + \int_0^T \int_\Omega \Big[ p(x,t) \, \partial_t q(x,t) 
      + \nabla_x p(x,t) \cdot \nabla_x q(x,t) \Big] dx \, dt \, .
\end{eqnarray*}
We end up with the variational problem to find $(u,p) \in X \times Y$ such
that
\begin{equation}\label{VF optimality system} 
  {\mathcal{B}}(u,p;v,q) =
  \langle u_T^\delta , q(T) \rangle_{L^2(\Omega)} \quad
  \forall \, (v,q) \in Y \times X,
\end{equation}
where the bilinear form ${\mathcal{B}}(\cdot,\cdot;\cdot,\cdot)$ is given as 
\[
  {\mathcal{B}}(u,p;v,q) := b(u,v) - b(q,p) +
  \langle u(T) , q(T) \rangle_{L^2(\Omega)} + \varrho \,
  \langle u(0) , q(0) \rangle_{L^2(\Omega)} \, .
\]
We note that the bilinear form $b(\cdot,\cdot)$, as defined by \eqref{Def
  b}, is bounded:  
\[
  |b(u,v)| \leq \sqrt{2} \, \| u \|_X \| v \|_Y \quad \forall \,
  u \in X, v \in Y.
\]
For $u \in X$ we have $\| u(0) \|_{L^2(\Omega)} \leq \mu \, \| u \|_X$ and
$\|u(T) \|_{L^2(\Omega)} \leq \mu \, \| u \|_X$ with 
\[
  \mu := \left(
    1 + \frac{1}{2} \left[ \frac{c_F}{T} \right]^2 +
    \sqrt{\frac{1}{4} \left[ \frac{c_F}{T} \right]^4 +
    \left[ \frac{c_F}{T} \right]^2}
    \right)^{1/2},
\]
where $c_F$ is the constant in Friedrichs' inequality in $H_0^1(\Omega)$. With
these ingredients, we are in the position to prove that the bilinear form
${\mathcal{B}}(\cdot,\cdot;\cdot,\cdot)$ is bounded, i.e., for all $(u,p),
(q,v) \in X \times Y$, there holds 
\[
  |{\mathcal{B}}(u,p;v,q)| \, \leq \, 2 \, (1+\varrho) \, \mu^2 \,
  \sqrt{\| u \|_X^2 + \| p \|_Y^2} \,
  \sqrt{\| q \|_X^2 + \| v \|_Y^2} \, .
\]
Moreover, we can establish the following inf-sup stability condition which can
be proved similarly to \cite[Lemma 3.2]{LangerOlafTroltzschYang01}. 
\begin{lemma}
  For simplicity, let us assume $\varrho  \in (0,1]$.
   Then there holds the inf-sup stability condition
  \[
    \frac{3}{10} \, \varrho \, \sqrt{\| u \|_X^2 + \| p \|_Y^2} \leq
    \sup\limits_{0 \neq (v,q) \in Y \times X}
    \frac{{\mathcal{B}}(u,p;v,q)}{\sqrt{\| q \|_X^2 + \| v \|_Y^2}}
    \quad \forall \, (u,p) \in X \times Y .
  \]
  Moreover, for all $(v,q) \in Y \times X$, there exist
  $(\overline{u}, \overline{p}) \in X \times Y$ satisfying
  \[
    {\mathcal{B}}(\overline{u},\overline{p};v,q) > 0 .
  \]
\end{lemma}

\noindent
Now, using the Banach--Ne\v{c}as--Babu\v{s}ka theorem
(see, e.g., \cite{ErnGuermond2004}), we can ensure well-posedness of
the variational optimality problem \eqref{VF optimality system} 
for any fixed positive regularization parameter $\varrho$.

%
%

\section{Space-time finite element methods}
\label{dist}

For the space-time finite element discretization of the
variational formulation \eqref{VF optimality system},
we first introduce conforming finite element spaces $X_h \subset X$
and $Y_h \subset Y$. In particular, we consider $X_h=Y_h$ 
spanned by piecewise linear continuous basis functions which are defined
with respect to some admissible decomposition of the space-time domain $Q$ into
shape regular simplicial finite elements. In addition, we will use the
subspace $Y_{0,h} \subset Y_h$ of basis functions with zero initial
values. Moreover, $Z_h \subset L^2(\Omega)$ is a finite element space to
discretize the control $z$. The space-time finite element discretization of
the forward problem \eqref{conheat} reads to find $u_h \in X_h$ such that 
\begin{equation}\label{VF forward}
  b(u_h,v_h) = 0 \quad \forall v_h \in Y_{0,h}, \quad
   \langle u_h - z_h , v_h \rangle_{L^2(\Sigma_0)} = 0 \quad
  \forall v_h \in Y_h \backslash Y_{0,h} .
\end{equation}
When denoting the degrees of freedom of $u_h$ at $\Sigma_0$, at $\Sigma_T$,
and in $Q$ by $\underline{u}_0$, $\underline{u}_T$, and $\underline{u}_I$,
respectively, the variational formulation 
\eqref{VF forward} is equivalent to the linear system
\[
  \left( \begin{array}{ccc}
           M_{00} & & \\ K_{0I} & K_{II} & K_{TI} \\ & K_{IT} & K_{TT}
  \end{array} \right)
  \left( \begin{array}{c}
  \underline{u}_0 \\ \underline{u}_I \\ \underline{u}_T \end{array} \right)
=
\left( \begin{array}{c} M_h^\top \underline{z} \\ \underline{0} \\
\underline{0} \end{array} \right),
\]
where the block entries of the stiffness matrix $K_h$ and the mass matrices
$M_{00}$ and $M_h$ are defined accordingly. After eliminating
$\underline{u}_0$, the resulting system corresponds to the space-time finite
element approach as considered in \cite{OS15}. In particular, we can
compute $\underline{u}_T = A_h \underline{z}$ to determine $u_h(T)$ in
dependency on the initial datum $z_h$, where  
\[
  A_h = \Big( K_{TT} - K_{IT} K_{II}^{-1} K_{TI} \Big)^{-1}
  K_{IT} K_{II}^{-1} K_{0I} M_{00}^{-1} M_h^\top =
  \widetilde{A}_h M_h^\top \, .
\]
Instead of the cost functional \eqref{fullopt}, we now
consider the discrete cost functional
\begin{eqnarray*}
  {\mathcal{J}}_{\varrho,h}(z_h)
  & = & \frac{1}{2} \, \| u_h(x,T) - u_T^\delta \|^2_{L^2(\Omega)} +
        \frac{\varrho}{2} \, \| z_h \|^2_{L^2(\Omega)} \\
  & = & \frac{1}{2} \, ( A_h^\top M_{TT} A_h \underline{z} , \underline{z} )
        - ( A_h^\top \underline{f} , \underline{z} ) +
        \frac{1}{2} \, \| u_T^\delta \|^2_{L^2(\Omega)} +
        \frac{\varrho}{2} \, ( \overline{M}_h \underline{z} , \underline{z}),
\end{eqnarray*}
whose minimizer is given as the solution of the linear system
\begin{equation}\label{ LGS z}
  A_h^\top (M_{TT} A_h \underline{z} - \underline{f}) +
  \varrho \, \overline{M}_h \underline{z} = \underline{0} .
\end{equation}
Note that $M_{TT}$ is the mass matrix formed by the basis functions of $X_h$
at $\Sigma_T$, $\overline{M}_h$ is the mass matrix related to the control
space $Z_h$, and $\underline{f}$ is the load vector of the target $u_T^\delta$
tested with basis functions from $X_h$ at $\Sigma_T$. When inserting
$\underline{u}_T = A_h \underline{z}$ and introducing
$\underline{p}_0 :=\widetilde{A}_h^\top (M_{TT} \underline{u}_T - \underline{f})$, 
$\underline{p}_T :=  ( K_{TT} - K_{IT} K_{II}^{-1} K_{TI} )^{-\top}( M_{TT} \underline{u}_T - \underline{f} )$,
$\underline{p}_I :=  - K_{II}^{-\top} K_{IT}^\top \underline{p}_T$,
this finally results in the linear system to be solved:
\begin{equation}\label{Linear system 1}
  \left(
    \begin{array}{ccccccc}
      & & & & - M_{00} & - K_{0I}^\top & \\
      & & & & & - K_{II}^\top & - K_{IT}^\top \\
      & & M_{TT} & & & - K_{TI}^\top & - K_{TT}^\top \\
      & & & \varrho \overline{M}_h & M_h & & \\
      M_{00} & & & - M_h^\top & & & \\
      K_{0I} & K_{II} & K_{TI} & & & & \\
       & K_{IT} & K_{TT} & & & &
    \end{array}
  \right)
  \left(
    \begin{array}{c}
      \underline{u}_0 \\
      \underline{u}_I \\
      \underline{u}_T \\
      \underline{z} \\
      \underline{p}_0 \\
      \underline{p}_I \\
      \underline{p}_T
    \end{array}
  \right)
  =
  \left(
    \begin{array}{c}
      \underline{0} \\
      \underline{0} \\
      \underline{f} \\
      \underline{0} \\
      \underline{0} \\
      \underline{0} \\
      \underline{0}
    \end{array}
  \right) .
\end{equation}
In the particular case, when $Z_h = Y_{h|\Sigma_0} \subset H^1_0(\Omega)$
is the space of piecewise linear basis functions as well, the mass
matrices $M_{00} = \overline{M}_h = M_h$ coincide, and therefore we
can eliminate $\underline{z} = \underline{u}_0$ and
$\underline{p}_0 = - \varrho \underline{z} = - \varrho \underline{u}_0$
to obtain
\begin{equation}\label{Linear system 2}
  \left( \begin{array}{ccccc}
  \varrho 
  M_{00} & & & - K_{0I}^\top & \\
  & & & - K_{II}^\top & - K_{IT}^\top \\     
  & & M_{TT} & - K_{TI}^\top & - K_{TT}^\top \\     
  K_{0I} & K_{II} & K_{TI} & & \\
  & K_{IT} & K_{TT} & & \\
  \end{array} \right)
  \left( \begin{array}{c}
           \underline{u}_0 \\ \underline{u}_I \\ \underline{u}_T  \\
           \underline{p}_I \\ \underline{p}_T 
  \end{array} \right) =
\left( \begin{array}{c}
         \underline{0} \\ \underline{0} \\ \underline{f} \\ \underline{0} \\
         \underline{0}
  \end{array} \right) .
\end{equation}
Note that \eqref{Linear system 2} is nothing but the Galerkin
discretization of the variational formulation
\eqref{VF optimality system} when using $X_h \subset X$ and $Y_{0,h}
\subset Y$ as finite element ansatz and test spaces. Obviously, the linear
system \eqref{ LGS z} and, therefore, \eqref{Linear system 2} are
uniquely solvable.  

In practice, the noise level $\delta \ge 0$ is usually given by the
measurement environment, and one has to choose suitable discretization and
regularization parameters $h$ and $\rho$. This is well investigated for linear
inverse problems; see, e.g., the classical book by Tikhonov and Arsenin
\cite{TikhonovArsenin:1977} and the more recent publications
\cite{EnglHankeNeubauer1996,GriesbaumKaltenbacherVexler2008}. In our
numerical experiments presented in the next section, we only play with the
parmeters $\delta$ and $h$ for a fixed small $\varrho$. 

%
%
\section{Numerical results}\label{numa}
\begin{figure}[ht]
   \centering
   \includegraphics[width=0.8\textwidth]{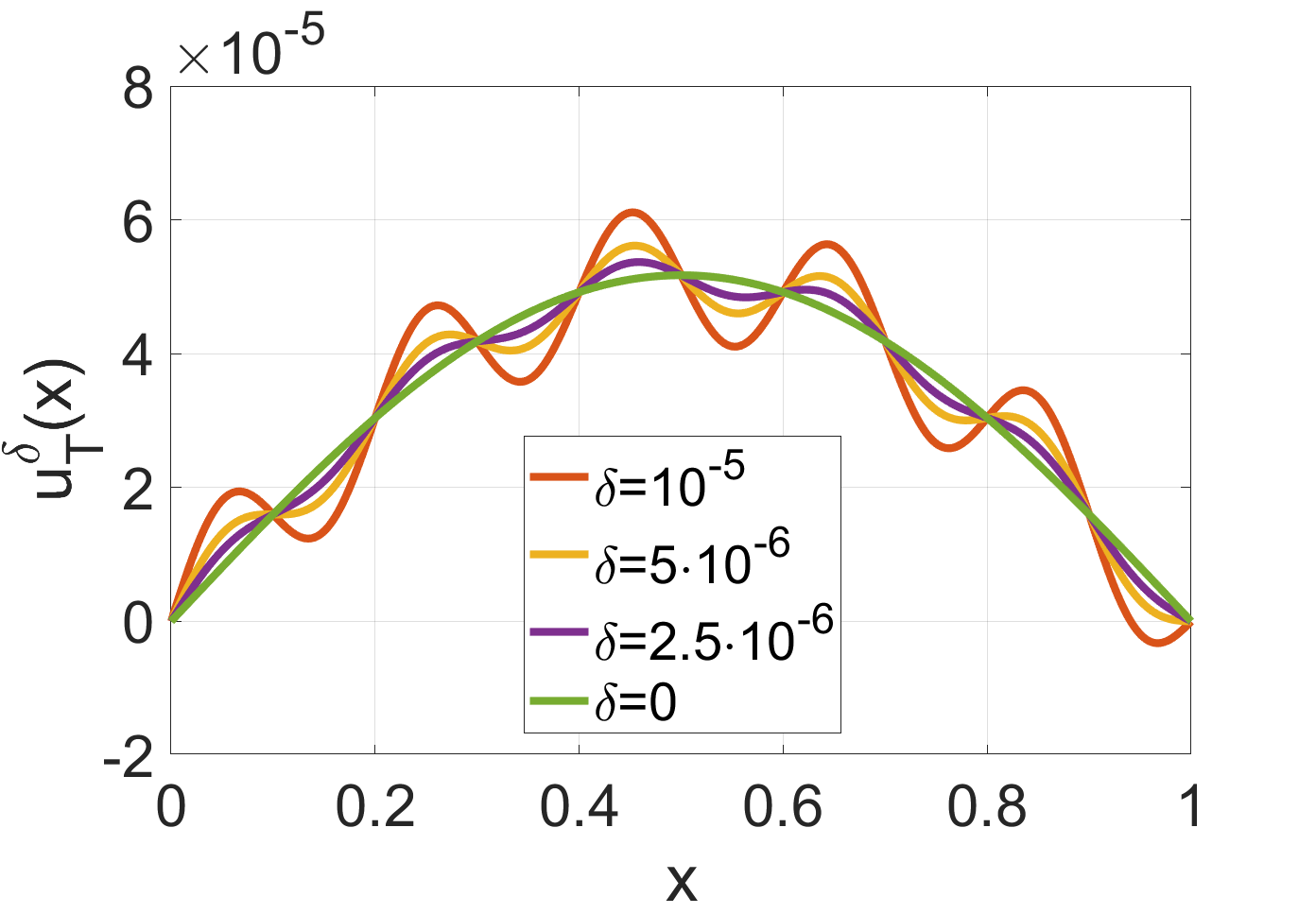}
   \caption{Comparison of the exact ($\delta=0$) and noisy
    ($\delta>0$) observation data.}\label{obsdata}
 \end{figure}

\begin{figure}[ht]
  \centering
  \includegraphics[width=0.48\textwidth]{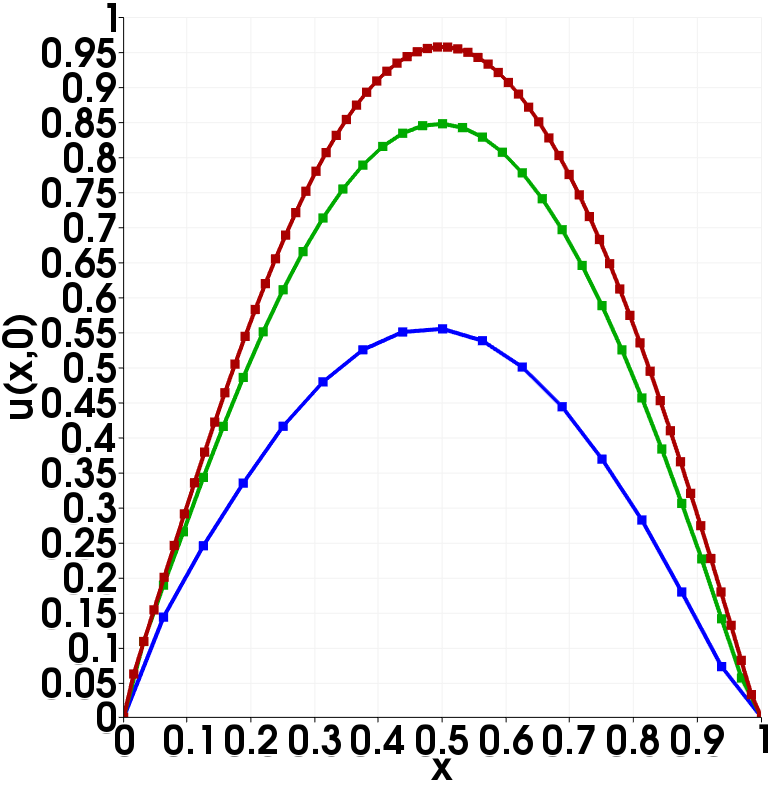}\,
  \includegraphics[width=0.48\textwidth]{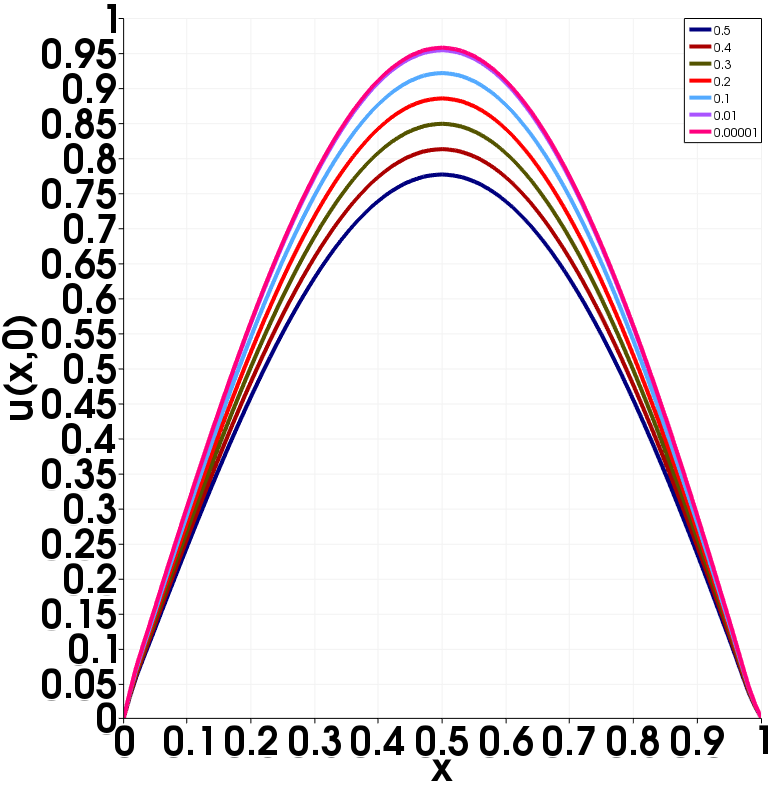}
  \caption{Convergence of the reconstructed initial data with respect
    to the mesh refinement $h\in\{1/16, 1/32, 1/64\}$, $\delta=0$,
    $\varrho=10^{-14}$ (left), and
    convergence with respect to the noise level
    $\delta\in\{0.5, 0.4, 0.3, 0.2, 10^{-1}, 10^{-3}, 10^{-5}\}$,
    $h=1/64$, $\varrho = 10^{-14}$ (right).}\label{recinitdelta0} 
\end{figure}

We take $\Omega=(0,1)$ and $T=1$, i.e., $Q=(0,1)^2$, and consider 
the manufactured observation data $u_T^\delta (x) := e^{-\pi^2}\sin(\pi
x)+\delta \sin(10\pi x)$ with some noise represented by the second term; see
exact and noisy data with $\delta\in\{0, 10^{-5}, 5\cdot 10^{-6}, 2.5\cdot
10^{-6}\}$ in Fig. \ref{obsdata}.  
To study the convergence of the space-time finite element solution to
the exact initial datum, we use the target $u_T(x)=e^{-\pi^2}\sin(\pi x)$
without any noise. The reconstructed initial data with respect to
a varying mesh size are illustrated in the left plot of
Fig. \ref{recinitdelta0}, where $\varrho=10^{-14}$. We clearly see the
convergence of the approximations to the exact initial datum with respect to
the mesh refinement. The right plot of Fig.~\ref{recinitdelta0} shows the
reconstructed initial approximation with different noise levels $\delta$. For
a decreasing $\delta$, we observe an improved reconstruction. 

%
%
\section{Conclusions}\label{sec:con}
We have applied the space-time FEM
from \cite{OS15} to the numerical solution of the classical inverse heat
conduction problem to determine the initial datum from measured observation
data at some time horizon $T$. The numerical results show the potential of
this approach for more interesting inverse problems. The space-time FEM are
very much suited for designing smart adaptive algorithms along the line
proposed in \cite{GriesbaumKaltenbacherVexler2008} determining the
optimal choice of $\varrho$ and $h$ for a given noise level $\delta$ in a
multilevel (nested iteration) setting.

\bibliography{Report}
\bibliographystyle{plain}

\end{document}